\documentclass{elsart}
\usepackage{amssymb}
\usepackage{latexsym}
\usepackage{amsmath}
\usepackage{eufrak}
\usepackage{theorem}

\font\sym=msbm10

\font\ssym=msbm8
\def\gD{\EuFfrak D}
\def\gH{\EuFrak H}

\def\gN{\EuFrak N}
\def\gD{\EuFrak D}

\def\mR{\hbox{\sym R}}
\def\smR{\hbox{\ssym R}}
\def\mC{\hbox{\sym C}}

\def\mZ{\hbox{\sym Z}}
\def\smo{\mathcal H}
\def\WF{M_{\smo,H}}
\newtheorem{Thr}{Theorem}
\newtheorem{Crl}{Corollary}
\newtheorem{Lm}{Lemma}
\newtheorem{Prp}{Proposition}
\begin{document}
\begin{frontmatter}
% Title, authors and addresses

% use the thanksref command within \title, \author or \address for footnotes;
% use the corauthref command within \author for corresponding author footnotes;
% use the ead command for the email address,
% and the form \ead[url] for the home page:
% \title{Title\thanksref{label1}}
% \thanks[label1]{}
% \author{Name\corauthref{cor1}\thanksref{label2}}
% \ead{email address}
% \ead[url]{home page}
% \thanks[label2]{}
% \corauth[cor1]{}
% \address{Address\thanksref{label3}}
% \thanks[label3]{}

\title{On  Periodic Matrix-Valued Weyl-Titchmarsh Functions}

% use optional labels to link authors explicitly to addresses:
% \author[label1,label2]{}
% \address[label1]{}
% \address[label2]{}

\author{M. Bekker}

\address{Department of Mathematics and Statistics,
University of Missouri-Rolla,\\ Rolla, MO 65409, USA}

\author{E.Tsekanovksii}

\address {Department of Mathematics, Niagara University,
 NY 14109, USA}
\begin{abstract}
% Text of abstract
We consider a certain class of Herglotz-Nevanlinna matrix-valued functions which 
can be realized as the Weyl-Titchmarsh
matrix-valued function of some symmetric operator and its self-adjoint extension.
New properties of Weyl -Titchmarsh matrix-valued functions as well as a new
version of the functional model in such realizations are presented. In the case of periodic
Herglotz-Nevanlinna matrix-valued functions we provide a complete characterization of their
realizations in terms of the corresponding functional model. We also obtain properties  
of a symmetric operator and its self-adjoint extension generating periodic Weyl-Titchmarsh
matrix-valued function.
We study pairs of operators (a symmetric operator and its self-adjoint
extension) with  constant Weyl-Titchmarsh matrix-valued 
functions and establish connections between such  pairs of 
operators and representations of the canonical commutation relations for 
unitary groups of operators in Weyl's form. As a consequence of such an approach 
we obtain the Stone-von Neumann theorem for two unitary groups of operators 
satisfying the commutation relations as well as some extension and refinement 
of the classical functional model for generators of those groups. Our examples 
include multiplication operators in weighted spaces, first and second order differential
operators, as well as the Schr\"odinger operator 
with  linear potential and its perturbation by
  bounded periodic potential.
\end{abstract}


\end{frontmatter}

% main text
%\section{}
%\label{}
\section {Introduction} In this paper we study a certain class of Herglotz-Nevanlinna 
matrix-valued functions which can be realized as the Weyl-Titchmarsh 
matrix-valued function $\WF(z)$ 
generated by the densely defined symmetric operator $\smo$ and its self-adjoint extension $H$
acting on some Hilbert space $\gH$ \cite{DM},\cite{GT}, \cite{GKMT}. The new properties of these functions as well as 
a new version of the functional model for the pair $(\smo,H)$ in terms of $\WF(z)$ 
are obtained. We introduce
 so called $(U,b)$-periodic pair of operators $(\smo,H)$, ( $U{\smo}U^*=\smo-bI$, 
  $UHU^*=H-bI$, $U$ is a unitary operator in $\gH$ ) and establish that the Weyl-Titchmarsh 
 matrix-valued function is $b$-periodic ( $\WF(z+b)=\WF(z)$ ) if and only if the corresponding
 pair of operators $(\smo,H)$ generating this matrix-valued function is $(U,b)$-periodic.
 It is shown that any Weyl-Titchmarsh function $\WF(z)$ corresponding to symmetric operator 
 $\smo$ with the defect indices $(1,1)$ which admits quasi-hermitian extension $\smo_v$ without 
 spectrum is always $\pi/tr(\Im\smo_v^{-1})$-periodic. Each $(U,b)-$periodic
 symmetric operator $\smo$ is associated with a group $\Gamma$ of transformations of the set
 $U(m)$ of all $m\times m$ unitary matrices into itself. It turned out that the
 group $\Gamma$ is cyclic if and only if an operator $\smo$ admits periodic extension.  
We consider pair of operators 
$(\smo,H)$ with the constant Weyl-Titchmarsh matrix-valued functions and find out connections
between such type of pairs and representations of the
canonical commutation relations for unitary 
groups of operators in Weyl's form. As a consequence of such approach we obtain the Stone-von Neumann 
theorem \cite{St} for two unitary groups of operators satisfying the
commutation relations as well as some extension
and refinement of the classical functional model for generators of those groups.
The examples of the Schr\"odinger operator with  linear potential and its perturbation by
bounded  periodic function and  are considered.
\section{The Weyl-Titchmarsh function.} Let $\gH$ be a Hilbert space, $\smo$
be a prime symmetric operator in $\gH$, that is $\gH$ does not contain a
proper subspace that reduces $\smo$, and in which $\smo$ induces self-adjoint
operator. Let  $\gD(\smo)$ denotes the domain of $\smo$. We assume that  
defect index of $\smo$ is $(m,m)$, $m<\infty$.  It means that for any non-real 
$z$ defect subspace $\gN_z=[(\smo -\bar{z}I)\gD(\smo)]^{\perp}$ has dimension 
$m$.   Let $H$ be a self-adjoint extension of $\smo$ in $\gH$ (an orthogonal 
extension) with domain $\gD(H)$. The Weyl-Titchmarsh function of the pair 
$(\smo,H)$, $\WF(z)$, is an operator-valued function whose values are operators
on $m-$dimensional space $\gN_i$. $\WF(z)$ is defined on the resolvent set
$\rho(H)$ of the operator $H$ by    
\begin{equation}
\WF(z)=P_+(zH+I)(H-zI)^{-1}|_{\gN_i}
\end{equation}
where $P_+$ is the orthogonal projection from $\gH$ onto $\gN_i$. 
From spectral representation of $H$ it follows that $\WF(z)$ can be written as
\begin{equation}\label{Wf}
M_{\mathcal{H},H}(z)=\int\limits_{\smR}\frac{\lambda z+1}{\lambda-z}d\sigma
(\lambda).
\end{equation}   
Values of a nondecreasing function  $\sigma(\lambda)$ are operators on $\gN_i$,
and it is defined by $\sigma(\lambda)=P_+E(\lambda)|_{\gN_i}$, where 
$E(\lambda)$ is the resolution of identity associated with $H$. We normalize 
$E(\lambda)$ by condition $E(\lambda)=1/2(E(\lambda+0)+E(\lambda-0))$. 
It is evidently that $\WF(z)$ is analytic on $\rho(H)$, particularly, for $\Im z\ne 0$,
and from (\ref{Wf} it follows that $\Im\WF(z)\ge 0$ for $z\in\mC_+$.
Therefore, $\WF(z)$ belongs to the
Heglotz-Nevalinna class. 

Function $\sigma$ has the following properties: 
 \begin{equation}\label{normalization1} 
\int_{\smR}d\sigma(\lambda)=I_{\gN_i};
\end{equation}
\begin{equation}\label{normalization2}
\int_{\smR}(1+\lambda^2)(d\sigma(\lambda)h,h)=\infty\qquad\forall h\in\gN_i,
\end{equation}
and $\sigma(\lambda)=1/2(\sigma(\lambda+0)+\sigma(\lambda-0))$. 
Condition (\ref{normalization1}) is obvious,  condition (\ref{normalization2}) 
follows from the fact, that according to von Neumann's formulas, for vector 
$h\in\gN_i$, $h\not\in\gD(H)$.  Condition (\ref{normalization1}) provides
normalization condition for the Weyl-Titvhmarsh function:  
 $\WF(i)=iI_{\gN_i}$.
From condition (\ref{normalization2}) it
follows that  points of growth of $\sigma$ form a noncompact set.

Selecting an orthonormal basis in $\gN_i$ we can identify the space $\gN_i$ 
with $\mC^m$, and regard $\WF(z)$ and $\sigma(\lambda)$ as operators on
$\mC^m$. Matrices of these operators with respect to the selected basis are
also denoted by $\WF(z)$ and $\sigma(\lambda)$. 
\vskip .4truecm 
Important property of the Weyl-Titchmatsh functions is given by the following
theorem.
\begin{Thr}. Let $\smo$ and $\tilde\smo$ be prime symmetric operators with
  equal defect numbers in Hilbert spaces $\gH$ and $\tilde\gH$ respectively, 
and $H$ and $\tilde H$ be their self-adjoint extensions. Suppose that there is 
the unitary operator $W:\gH\rightarrow\tilde\gH$ such that 
$W\smo=\tilde\smo W$ and $WH=\tilde HW$. Then there is a unitary operator
$W_0:\gN_i\rightarrow\tilde\gN_i$ such that $W_0\WF(z)=M_{\tilde\smo,\tilde
  H}(z)W_0$.   
\end{Thr}
{\bf Proof.} From the assumptions of the Theorem it follows that $WE(\lambda)=
\tilde E(\lambda)W$, where $E(\lambda)$ and $\tilde E(\lambda)$ are the
resolutions of identity, associated with $H$ and $\tilde H$ respectively. It
is also obvious that $W\gN_z=\tilde\gN_z$, and $WP_+=\tilde P_+W$. 

Put $W_0=W|\gN_i$. Then $W_0$ is the unitary operator from $\gN_i$ onto 
$\tilde\gN_i$, $W_0^*=W^*|\tilde\gN_i$. For any $f\in\gN_i$ and $\tilde
g\in\tilde\gN_i$ we have
\begin{multline*}
(W_0\WF(z)f,\tilde g)=(W\WF(z)f,\tilde g)=(\WF(z)f,W^*\tilde g)=\\
\int\limits_{\mR}\frac{\lambda z+1}{\lambda-z}d(P_+E(\lambda)f,W^*\tilde g)=
\int\limits_{\mR}\frac{\lambda z+1}{\lambda-z}d(WP_+E(\lambda)f,\tilde g)=\\
\int\limits_{\mR}\frac{\lambda z+1}{\lambda-z}d(\tilde P_+\tilde E(\lambda)Wf,
\tilde g)= (M_{\tilde\smo,\tilde H}(z)W_0f,\tilde g). 
\end{multline*}
These equalities show that $W_0$ possesses desired property. 

If $\{e_j\}_{j=1}^m$ is an arbitrary orthonormal basis in $\gN_i$, then 
$\{W_0e_j\}$ is the orthonormal basis in $\tilde\gN_i.$ With respect to these
bases matrices of $\WF(z)$ and  $M_{\tilde\smo,\tilde H}(z)$ are equal.
Therefore, the Theorem 1 can be reformulated as following:

\noindent{\it If pairs $(\smo, H)$ and $(\tilde\smo,\tilde H)$ are unitarily 
equivalent, then there are bases with respect to which matrices of their Weyl-
Titchmarsh functions are equal}.      
%\vskip .5truecm
The next Theorem is the statement about realization. It provides the
functional model of the pair with prescribed Weyl-Titchmarsh function.
\begin{Thr} Let $F(z)$ be a function whose values are linear operators on
  $m-$\\dimensional space $\gN$, and  which admits integral representation 
\begin{equation*}
F(z)=\int\limits_{-\infty}^{\infty}\frac{\lambda
  z+1}{\lambda-z}d\sigma(\lambda)
\end{equation*}
where $\sigma(\lambda)$ is a nondecreasing function with values on the 
set of linear operators on $\gN$, and which satisfies
(\ref{normalization1})
and (\ref{normalization2}). Then there 
are Hilbert space $\tilde\gH$, prime symmetric operator $\tilde\smo$ with defect index 
$(m,m)$, and its self-adjoint extension $\tilde H$ in $\tilde\gH$, such that 
$F(z)=M_{\mathcal{\tilde H},\tilde H}(z)$. 
If $(\hat\gH,\hat\smo,\hat H)$ is another realization of $F$, then there
is a unitary operator $\Psi:\tilde\gH\to\hat\gH$ such that $\Psi\tilde\smo=\hat\smo \Psi$, and
$\Psi\tilde H=\hat H\Psi$.
\end{Thr}
{\bf Proof.} Since $\sigma(\lambda)$ is nondecreasing operator-function 
 and satisfies (\ref{normalization1}), it is the generalized resolution of
identity which acts in $\gN$. According to the theorem of M.A.Najmark (see, 
for example \cite{AkG}) there exist a Hilbert space $\tilde\gH$ which contains 
$\gN$ as a subspace and the orthogonal resolution of identity $\tilde E(\lambda)$, 
such that for any Borel set
$\Delta\in\mathcal B(\mR)$ ($\mathcal B(\mR)$ is the Borel field of $\mR$)
$\sigma(\Delta)=P\tilde E(\Delta)|\gN$, where $P$ is the orthogonal projection from
$\tilde\gH$ onto $\gN$. The space $\tilde\gH$ can be selected minimal in that sense that
$c.l.h.\{\tilde E(\Delta)h|\Delta\in\mathcal B(\mR),h\in\gN\}=\tilde\gH$, where 
{\it  c.l.h} means closed linear hall. The orthogonal resolution of identity
$\tilde E(\lambda)$ defines the  self-adjoint operator $H$ in $\tilde\gH$. Under minimality
condition the Hilbert space $\tilde\gH$ and the operator $\tilde H$ are defined uniquely up 
to unitary equivalence. 

In our situation this construction gives  Hilbert space
$\tilde\gH=L^2(\mR,\gN,d\sigma)$. 
Elements of $\tilde\gH$ are measurable functions $f(\lambda)$, $\lambda\in\mR$ 
with values in $\gN$ such that 
\begin{equation*}
\int\limits_{\smR}(d\sigma(\lambda)f(\lambda),f(\lambda))_{\gN}<\infty.
\end{equation*}
The space $\gN$ is identified with subspace of $L^2(\mR,\gN,d\sigma)$ which
consists of constant functions.
The orthogonal resolution of identity $\tilde E$ is defined as 
$\tilde E(\Delta)f(\lambda)=\chi_{\Delta}(\lambda)f(\lambda)$, where $\chi_{\Delta}$
is the indicator function of the set $\Delta$.

 The self-adjoint operator $\tilde H$ defined as follows:
\begin{equation}\label{domself}
\gD(\tilde H)=\{f\in\tilde\gH|\int_{\smR}(1+\lambda^2)(d\sigma(\lambda)f(\lambda),f(\lambda))_{\gN}<\infty\},
\end{equation}
\begin{equation}\label{self}
(\tilde Hf)(\lambda)=\lambda f(\lambda),\qquad f\in\gD(\tilde H).
\end{equation}
From (\ref{normalization2}) it follows that $\tilde H$ is unbounded operator.

Put
\begin{equation} \label{ds}
\gD(\tilde\smo)=\{f\in\gD(\tilde H)|\int_{\smR}(\lambda+i)d\sigma(\lambda)f(\lambda)=0\},
\end{equation}
and 
\begin{equation}
(\tilde\smo f)(\lambda)=\lambda f(\lambda),\qquad\lambda\in\gD(\smo).
\end{equation}
$\gD(\tilde\smo)$ is linear manifold, dense in $\tilde\gH$ (this fact follows from
(\ref{normalization2})), 
and $(\tilde\smo f,g)=(f,\tilde\smo g)$ for $f,g\in\gD(\tilde\smo)$. Thus, $\tilde\smo$ is a 
symmetric operator. Moreover, 
condition (\ref{ds}) implies, that 
$\gN=[(\tilde\smo+iI)\gD(\tilde\smo)]^{\perp}=\gN_i$. Indeed, for $f\in
L^2(\mR,\gN,d\sigma)$ put $f_0=\int d\sigma(\lambda)f$. Then we 
 have $f=(\lambda+i)g+h$, where $g=(f-f_0)/(\lambda+i)\in\gD(\tilde\smo)$, 
$h=f_0\perp (\lambda+i)g$. Therefore, one of the defect 
numbers of $\tilde\smo$ is $m$. It is easily seen, that 
$\gN_{-i}=\{\displaystyle{\frac{\lambda-i}{\lambda+i}}\xi|\xi\in\gN\}$, 
which means that $\dim\gN_{-i}=m$, and defect index of $\tilde\smo$ is $(m,m)$. 
In general, for arbitrary nonreal $z$ the defect subspace 
$\gN_z=\{\displaystyle{\frac{\lambda-i}{\lambda-z}}\xi|\xi\in\gN\}$.

The Weyl-Titchmarsh function for the pair $(\tilde\smo,\tilde H)$ is 
\begin{equation*}
M_{\tilde\smo,\tilde H}=P_+(zH+I)(H-zI)^{-1}|_{\gN_i}=\int\limits_{\smR}\frac{z\lambda+1}
{\lambda-z}d\sigma(\lambda)
\end{equation*}
coincides with the given function $F$.  Uniqueness of this realization (up to
unitary equivalence) is provided by the Najmark's theorem.

Combining results of the Theorems 1 and 2 we obtain the following statement
(see \cite{GT}, \cite {GKMT}).
\begin{Crl}. Let $\smo$ be a prime symmetric operator on a Hilbert space $\gH$
with index of defect $(m,m)$ ($m<\infty$), and let $H$ be a self-adjoint
extension of $\smo$ in $\gH$. Let $\WF(z)$ be the Weyl-Titchmarsh function
of the pair $(\smo,H)$. Let $(\tilde\gH,\tilde\smo,\tilde H)$ be the
realization of $\WF$ described in Theorem 2. Then there is the unitary
operator $\Phi:\gH\to\tilde\gH$ such that  
\begin{equation}
\tilde\smo=\Phi\smo\Phi^*,
\end{equation}
and
\begin{equation}
\tilde H=\Phi H\Phi^*.
\end{equation}
\end{Crl}
Let $U$ be a unitary operator on $\gH$, and $\tilde U=\Phi U\Phi^*$ be its
representation in the model space $\tilde\gH$. We say that the operator $U$ is 
of {\bf shift-type} (s-type) operator if for $f\in\tilde\gH$
\begin{equation}\label{stype}
(\tilde Uf)(\lambda)=D\frac{\lambda-i}{\lambda-i-b}f(\lambda-b),
\end{equation}
where $D$ is a unitary operator on $\gN$ which commutes with
$\sigma(\lambda)$,  and $b$ is a real number.

Often it is more convenient to use the following realization of $F$ (see \cite
 {GT},\cite{GKMT}). 
Put 
\begin{equation}\label{tau}
d\tau(\lambda)=(1+\lambda^2)d\sigma(\lambda).
\end{equation} 
Then 
\begin{equation}\label{rfunction}
F(z)=\int\limits_{-\infty}^{\infty}\left[\frac{1}{\lambda-z}-\frac{\lambda}
{1+\lambda^2}\right]d\tau(\lambda)
\end{equation} 

The mapping $W:L^2(\mR,\gN,d\sigma)\to L^2(\mR,\gN,d\tau)$,
$(Wf)(\lambda)=f(\lambda)/(\lambda-i)$
is a unitary one . For the self-adjoint operator $\hat H=WHW^*$ we have

\begin{equation*}
\gD(\hat H)=\{f\in L^2(\mR,\gN,d\tau)|\int_{\smR}(1+\lambda^2)(d\tau
(\lambda)f(\lambda),f(\lambda))_{\gN}<\infty\},
\end{equation*}
and $\hat Hf(\lambda)=\lambda f(\lambda)$.
 
For symmetric operator $\hat\smo=W\smo W^*$  

\noindent ({\it i}) $\gD(\hat\smo)=\{f\in\gD(\hat H)|\int_{\smR}f(\lambda)d\tau(\lambda)=0$;
\\
\noindent ({\it ii}) $(\hat\smo f)(\lambda)=\lambda f(\lambda)$;
\\
\noindent $\gN_z=\{\displaystyle{\frac{1}{\lambda-z}\xi}|\xi\in\gN\}$.\\

In such representation the s-type unitary operator $U$ acts as \\$(\hat
Uf)(\lambda)=Df(\lambda-b)$.

Some additional properties of  the Weyl-Titchmarsh functions 
and their applications can be found in \cite{GT}, \cite{GKMT}.

\section{Periodic Operators.}
Let $\smo$ be a prime symmetric operator with index of defect $(m,m)$, 
$m<\infty$ and $H$ be its orthogonal self-adjoint extension. In this section we 
study pairs $(\smo,H)$ for which the Weyl-Titchmarsh function is $b-$periodic, 
that is 
\begin{equation}\label{bperiodic}
\WF(z)=\WF(z+b),
\end{equation} 
where $b$ is some real number.  
 
We start from the following  lemma.
\begin{Lm}. Let $F(z)$ be a function whose values are linear operators on
  $m-$\\dimensional space $\gN$, and  which admits integral representation 
\begin{eqnarray*}
F(z)=\int\limits_{-\infty}^{\infty}\frac{\lambda
  z+1}{\lambda-z}d\sigma(\lambda)=zI_{\gN}+(1+z^2)
\int\limits_{-\infty}^{\infty}\frac{1}{\lambda-z}d\sigma(\lambda),
\end{eqnarray*}
where $\sigma(\lambda)$ is a nondecreasing function with values on 
the set of linear operators on $\gN$ which satisfies conditions 
(\ref{normalization1}) and (\ref{normalization2}). 
The function $F(z)$ is $b-$periodic, if and only
 if
\begin{equation} \label{tper} 
\tau(\Delta+b)=\tau(\Delta) 
\end{equation}
for any 
$\Delta\in\mathcal {B}(\mR)$, where 
 $\tau$ is defined by (\ref{tau}).
\end{Lm}
{\bf Proof.} In order to prove the Lemma we need the following generalization 
of the Stieltjes inversion formula. This generalization due to M.Livsic 
(see \cite{Kr}, Lemma 2.1):
%\vskip .3truecm

{\it Let $\sigma(\lambda)=1/2(\sigma(\lambda+0)+\sigma(\lambda-0))$ 
$(-\infty<\lambda<\infty)$ be some function of bounded variation on each
finite interval, such that the integral
\begin{equation*}
\Phi(z)=\int\limits_{-\infty}^{\infty}\frac{d\sigma(\lambda)}{\lambda -z}
\end{equation*}
converges absolutely.

Let $\varphi(\lambda)$ be some function analytic on the closed interval
$\Delta=[\alpha,\beta]$.

Denote by $\Delta_{\epsilon}$ the broken path of integration consisting of
directed segment $[\alpha-i\epsilon,\beta-i\epsilon]$ and antiparallel segment 
$[\beta+i\epsilon, \alpha+i\epsilon]$.

Then
\begin{equation*}
\lim\limits_{\epsilon\to 0}\frac{1}{2\pi
  i}\int\limits_{\Delta_{\epsilon}}\varphi(z)\Phi(z)dz=
-\int\limits_{\alpha}^{\beta}\varphi(\lambda)d\sigma(\lambda).
\end{equation*} 
}

Fix and orthonormal basis $\{e_j\}_{j=1}^m$ in the space $\gN$.
Condition of $b-$periodicity of the function $F(z)$
gives
\begin{equation}\label{stieltjes}
b\delta_{jk}+(1+(z+b)^2)\int\limits_{-\infty}^{\infty}\frac{1}{\lambda-b-z}d\sigma_{jk}(\lambda)=
(1+z^2)\int\limits_{-\infty}^{\infty}\frac{1}{\lambda-z}d\sigma_{jk}(\lambda),
\end{equation}
$\Im z\ne 0$, and $\sigma_{jk}(\lambda)=(\sigma(\lambda)e_k,e_j)$. Since
$dim\gN=m<\infty$, variations of all functions $\sigma_{jk}$,
$j,k=1,2,\ldots,m$ are uniformly bounded and (\ref{tper}) follows from the
Livsic's lemma. Indeed, evaluating the integral of both sides of
(\ref{stieltjes}) along $\Delta_{\epsilon}$ and then taking the limit as
$\epsilon\to 0$ we obtain
\begin{equation*}
\int\limits_{\alpha}^{\beta}[1+(\lambda+b)^2]d\sigma(\lambda+b)=\int\limits_{\alpha}^{\beta}(1+\lambda^2)d
\sigma(\lambda),
\end{equation*} 
which is (\ref{tper}).

Suppose now that (\ref{tper}) is fulfilled. Then we have
\begin{eqnarray*}
F(z+b)-F(z)=\int\limits_{\smR}\left[\frac{1}{\lambda-z-b}-\frac{1}{\lambda-z}\right] 
d\tau(\lambda)=c,
\end{eqnarray*} 
$c=\int_{\smR}[\lambda/(1+\lambda^2)-(\lambda+b)/(1+(\lambda+b)^2)]d\tau(\lambda)$,
and the integrals converge absolutely.  We assume for simplicity that
$m=1$ (for case $m<\infty$ the  proof can be done by componentwise
arguments). Consider difference 
\begin{eqnarray*}
|F(iy+b)-F(iy)|=\left|\int\limits_{\smR}\left[\frac{1}{\lambda-iy-b}-
\frac{1}{\lambda-iy}\right]d\tau(\lambda)\right|\le\\
 b\int\limits_{\smR}\frac{d\tau(\lambda)}{\sqrt{\lambda^2+y^2}
\sqrt{(\lambda-b)^2+y^2}}.   
\end{eqnarray*}
For large $y$ we have $1/(\sqrt{\lambda^2+y^2}\sqrt{(\lambda-b)^2+y^2})\le 
1/(\sqrt{\lambda^2+1}\sqrt{(\lambda-b)^2+1})$, therefore, there is $A>0$ such
that 
\begin{equation*}
\left(\int\limits_{-\infty}^{-A}+\int\limits_{A}^{\infty}\right)\frac{d\tau(\lambda)}{\sqrt{\lambda^2+y^2}
\sqrt{(\lambda-b)^2+y^2}}<\frac{\epsilon}{2}
\end{equation*}
for any $\epsilon >0$ uniformly with respect to $y$. From the other side,
\begin{eqnarray*}
\int\limits_{-A}^{A}\frac{d\tau(\lambda)}{\sqrt{\lambda^2+y^2}
\sqrt{(\lambda-b)^2+y^2}}\le\frac{1+A^2}{y^2}<\frac{\epsilon}{2}
\end{eqnarray*} 
for $y$ large enough (we used the fact that
$\int_{\smR}d\tau(\lambda)/(1+\lambda^2)=1$). We have proved that the
constant $c=0$, and $F(z+b)=F(z)$.

The Lemma is proved now.
%\vskip .5truecm

{\bf Definition.} An operator $T$ acting on a Hilbert space $\gH$  with domain
$\gD(T)$ is said to be  $(U,b)-${\bf periodic}, if there  is a unitary operator $U$ such that 
\begin{equation}\label{periodic1}
U\gD(T)\subset\gD(T),
\end{equation}
\begin{equation}\label{periodic2}
UTU^*=T-bI
\end{equation}
for some number $b$. 

Of course, periodic operator cannot be bounded. One can
easily see that if the operator $T^*$ exists, then it is $(U,\bar b)-$periodic.
\vskip .3truecm
We say that prime symmetric operator $\smo$ in $\gH$ and its self-adjoint 
extension $H$ form a $(U,b)$-{\bf periodic pair}, if conditions 
(\ref{periodic1}) and (\ref{periodic2}) are  fulfilled for both of them 
(with the same unitary operator $U$). 

It is evidently, that if $\smo$ is a $(U,b)-$periodic periodic operator,
and $\gN_z$ is a defect subspace of $\smo$, then $U\gN_z=\gN_{z+b}$. 
\begin{Prp}. Let $\smo$ be a prime symmetric operator, $H\supset\smo$ be its
  selfadjoint extension such that the pair 
$(\smo, H)$ is $(U,b)-$periodic and $(V,b)$-periodic. Then the unitary
operator $W=V^*U$ has following properties:
\begin{enumerate}
\item
$W$ commutes with $H$;
\item  each defect subspace $\gN_z$ reduces W;
\item if $\smo$ has  defect index $(m,m)$, $m<\infty$, then the spectrum of $W$
consists of finite number of eigenvalues; number of distinct eigenvalues not
greater than $m$. 
\end{enumerate} 
\end{Prp}
Indeed, properties {\it 1} and {\it 2} follow directly from the definitions
above. The property {\it 3} follows from the fact that the operator $W$
commutes with the resolution of identity $E(\lambda)$ associated with $H$, 
$c.l.h.\{E(\Delta)\gN|\Delta\in{\mathcal B}(\mR)\}=\gH$, where $\gN$ is a
defect subspace of $\smo$, and the spectrum of $W|\gN$ consists of finite
numebrs of eigenvalues.

\begin{Thr} Let $\smo$ be a prime symmetric operator on a Hilbert space $\gH$ 
with defect index $(m,m)$, ($m<\infty$), and let $H$ be its self-adjoint 
extension in $\gH$. Then the following conditions are equivalent:
\begin{enumerate}
\item
The Weyl-Titchmarsh function $\WF(z)$ of the pair $(\smo,H)$ is $b-$periodic;
\item
The pair $(\smo,H)$ is $(U,b)-$periodic, where $U$ is an s-type operator.
\end{enumerate}
\end{Thr}
{\bf Proof.} Let pair $(\smo,H)$ has $b-$periodic Weyl-Titchmarsh function. 
Let $(\tilde\gH,\tilde\smo,\tilde H)$ be the realization of  
$(\gH,\smo,H)$, described in the Theorem  2. According to the Lemma 1  the
function $\sigma(\lambda)$ satisfies 
the periodicity condition \\$(1+(\lambda+b)^2)d\sigma(\lambda+b)=(1+\lambda^2)
d\sigma(\lambda)$. On the space $\tilde\gH=L^2(\mR,\gN_i,d\sigma)$ consider the 
operator $\tilde U:f\to \tilde Uf$ defined by
\begin{equation}
(\tilde Uf)(\lambda)=\frac{\lambda-i}{\lambda-b-i}f(\lambda-b).
\end{equation}
Operator $\tilde U$ is a unitary operator in $L^2(\mR,\gN_i,d\sigma)$. Indeed,  
\begin{eqnarray*}
(\tilde Uf,\tilde Uf)=\int\limits_{-\infty}^{\infty}\frac{\lambda^2+1}{1+(\lambda-b)^2}
(d\sigma(\lambda)f(\lambda-b),f(\lambda-b))=\\
\int\limits_{-\infty}^{\infty}\frac{1+(\lambda-b)^2}{1+(\lambda-b)^2}
d(\sigma(\lambda-b)f(\lambda-b),f(\lambda-b))=(f,f).
\end{eqnarray*}
The domain of the operator $\tilde\smo$ is invariant under $\tilde U$. For
$f\in\gD(\tilde \smo)$, that is $\int_{\smR}(\lambda+i)d\sigma(\lambda)f(\lambda)=0$, we have
\begin{eqnarray*}
\int\limits_{-\infty}^{\infty}(\lambda+i)d\sigma(\lambda)(Uf)(\lambda)=
\int\limits_{-\infty}^{\infty}\frac{\lambda^2+1}{\lambda-i-b}d\sigma(\lambda)
f(\lambda-b)=\\
\int\limits_{-\infty}^{\infty}\frac{1+(\lambda-b)^2}{\lambda-b-i}
d\sigma(\lambda-b)f(\lambda-b)=\int\limits_{-\infty}^{\infty}(\lambda+i)d\sigma(\lambda)f(\lambda)=0.
\end{eqnarray*}
It is obvious that if $f\in\gD(\tilde H)$, then $\tilde U\tilde Hf=(\tilde
H-bI)\tilde Uf$. Therefore, $(\tilde\smo,\tilde H)$ 
is the $(\tilde U,b)-$periodic pair. Therefore, the pair $(\smo, H)$ is the
$(U,b)-$periodic one, and $U$ is the s-type operator.  

Conversely, let $(\smo, H)$ be a $(U,b)$-periodic pair, with  operator $U$ 
of s-type. Therefore, in the realization $(\tilde\gH,\tilde\smo,\tilde H)$ the
pair $(\tilde\smo,\tilde H)$ is $(\tilde U,b)-$periodic, with $\tilde U$ of
the form (\ref{stype}).
From the equation 
$\tilde U\tilde H\tilde U^*=\tilde H-bI$ it follows that the resolution of
identity $\tilde E(\lambda)$ of the operator $\tilde H$ satisfies the condition 
\begin{equation}\label{spectralfunction}
\tilde U\tilde E(\lambda)\tilde U^*=\tilde E(\lambda+b).
\end{equation}
If $\hat\gN_i$ is the defect subspace of the operator $\tilde U\tilde H\tilde
U^*$, then $\hat\gN_i=\gN_{i+b}$. Let $\{e_j\}$ be an orthonormal basis in
$\gN$. Then $\tilde Ue_j=\displaystyle{\frac{\lambda-i}{\lambda-i-b}}De_j$, $j=1,2\ldots,m$
is the orthonormal basis in $\hat\gN_i=\gN_{i+b}$. Now the Theorem 1 gives
\begin{eqnarray*}
\sigma_{jk}(\lambda)=(\tilde E(\lambda)e_k,e_j)=(\tilde E(\lambda+b)\tilde
Ue_k,\tilde Ue_j)=\\
\int\limits_{-\infty}^{\lambda+b}\frac{1+s^2}{1+(s-b)^2}d\sigma(s),
\end{eqnarray*}
from which we get  $(1+\lambda^2)d\sigma(\lambda)=(1+(\lambda+b)^2)d\sigma(\lambda+b)$.

Therefore, the function $\sigma$ satisfies the condition of the Lemma 1, and 
$\WF(z)$ is the $b-$periodic function.
The Theorem is proved.  

{\it Remark}. It can be proved, that if $(\smo,H)$ is a $(U,b)-$periodic pair,
where index of defect of $\smo$ is $(1,1)$, then the unitary operator $U$ is
necessarily of s-type.

\begin{Lm}. Let $\smo$ be a $(U,b)$-periodic prime symmetric operator with 
finite and equal defect numbers, and let $(\smo, H_0)$ is a $(U,b)$-periodic
pair. Define operator functions $\mathcal{A}(z)$ and $\mathcal{B}(z)$ by the equations 
\begin{equation}\label{Az}
{\mathcal A}(z)=\int_{\smR}\frac{\lambda -i}{\lambda -z}d\sigma_0(\lambda),
\end{equation}
 \begin{equation}\label{Bz}
{\mathcal B}(z)=\int_{\smR}\frac{\lambda +i}{\lambda -z}d\sigma_0(\lambda),
\end{equation}
where $\sigma_0(\lambda)=P_+E_0(\lambda)|_{\gN_i}$, $E_0(\lambda)$ is the 
resolution of identity for $H_0$. Then the functions $\mathcal A$ and 
$\mathcal B$ satisfy the following identities:
\begin{equation}\label{Aztransform}
{\mathcal A}(z+b)=\frac{z+i}{z+b+i}{\mathcal A}(z),
\end{equation}
\begin{equation}\label{Bztransform}
{\mathcal B}(z+b)=\frac{z-i}{z+b-i}{\mathcal B}(z).
\end{equation}

\end{Lm}
{\bf Proof.} We prove identity for $\mathcal A$. Identity for $\mathcal B$ is proved similarly.

\begin{eqnarray*}
&{\mathcal A}(z+b)=
\int\frac{\lambda-i}{\lambda-z-b}d\sigma_0(\lambda)=&\\
&\frac{1}{z+b+i}\int\left[\frac{1}{\lambda-z-b}-\frac{1}{\lambda+i}\right]
(1+\lambda^2)d\sigma_0(\lambda).&
\end{eqnarray*} 
Since $(\smo, H_0)$ is the $(U,b)$-periodic pair, the Weyl-Titchmarsh function\\ 
$M_{\smo, H_0}(z)$ for the pair has period $b$, from which it follows, that
the measure $d\tau_0(\lambda)=(1+\lambda^2)d\sigma_0(\lambda)$ also has period 
$b$. This condition provides, that
\begin{equation*}
\int\left[\frac{1}{\lambda-z-b}-\frac{1}{\lambda+i}\right]d\tau_0(\lambda)
=\int\left[\frac{1}{\lambda-z}-\frac{1}{\lambda+i}\right]d\tau_0(\lambda),
\end{equation*} 
and the statement regarding the function $\mathcal A(z)$ follows. 
 
\begin{Crl}. Let $\smo$ be a prime symmetric operator in the Hilbert space 
$\gH$  with index of defect $(m,m)$, and $H_0$ be its orthogonal self-adjoint 
extension such that the pair $(\smo, H_0)$ is a $(U,b)$-periodic. Then for any 
other orthogonal self-adjoint extension $H$ of the operator $\smo$
corresponding  pair $(\smo, H)$ is  a $(U^{\prime},b)$-periodic with some  
unitary operator $U^{\prime}$.
\end{Crl}  
{\bf Proof.} In light of the Theorem 1 it is enough to show that periodicity of $M_{\smo, H_0}(z)$ implies periodicity of $\WF(z)$. 

Let $\sigma_0$ be the non decreasing operator valued function which provides 
 the integral representation of the  $M_{\smo, H_0}(z)$. 
Consider the functional model for the pair $(\smo, H_0)$. 
Then the domain $\gD(H)$ of the self-adjoint extension $H$ of the operator
 $\smo$ consists of the functions $f(\lambda)\in L^2(\mR,\gN_i,d\sigma_0)$ which 
can be written as 
\begin{equation}
f=g+(\varphi_i-V\varphi_{-i}),
\end{equation}
where $g\in\gD(\smo)$, that is $\int_{\smR}(\lambda+i)g(\lambda)d\sigma_0(\lambda)=0$, $\varphi\in\gN_i$, $\varphi_{-i}\in\gN_{-i}$, 
$\Vert\varphi\Vert=\Vert\varphi_{-i}\Vert$, and $V$ is a some unitary operator
 in $\gN_{-i}$. We also have that for $f\in\gD(H)$ $Hf=\smo g+i(\varphi_i+V\varphi_{-i})$. 

From the definition of Weyl-Titchmarsh function of the pair we have that 
\begin{equation*}
\frac{\WF(z)-M_{\smo, H_0}(z)}{1+z^2}=P_+\left[R(z)-R_0(z)\right]|_{\gN_i},
\end{equation*}
where $R$ and $R_0$ are resolvents of $H$ and $H_0$ respectively. Calculating 
the difference of resolvents, we get  the following expression
\begin{equation}
\frac{\WF(z)-M_{\smo, H_0}(z)}{1+z^2}={\mathcal A}(z)\left(I-V\right)
\left[(i+z){\mathcal A}(z)V+(i-z){\mathcal B}(z)\right]^{-1}{\mathcal B}(z),
\end{equation}
where ${\mathcal A}(z)$ and ${\mathcal B}(z)$ are defined by (\ref{Az}) and 
(\ref{Bz}). Using now formulas (\ref{Aztransform}) and (\ref{Bztransform}), 
we obtain that  $\WF(z)-M_{\smo, H_0}(z)=\WF(z+b)-M_{\smo, H_0}(z+b)$, and the Corollary is proved.
\vskip 1.0truecm
Let $\smo$ be a $(U,b)-$periodic prime symmetric operator in a Hilbert space 
$\gH$ with index of defect $(m,m), (m<\infty)$. Fix orthonormal bases 
$\{\varphi_j\}_{j=1}^m$ in $\gN_i$ and $\{\psi_j\}_{j=1}^m$ in $\gN_{-i}$, 
and a unitary operator $V_0$ in $\gN_{-i}$. The matrix of this operator with 
respect to the basis $\{\psi_j\}_{j=1}^m$ we also denote by $V_0$. Denote by 
$\gD(H_0)$ the domain of self-adjoint extension $H_0$ of the operator $\smo$ 
defined as 
\begin{equation*}
\gD(H_0)=\{f\in\gH|f=f_0+\sum_{j}c_j(\varphi_j-V_0\psi_j), f_0\in\gD(\smo),
c_j\in\mC\}.
\end{equation*}
Since $U\gD(\smo)=\gD(\smo)$ the set $U^n\gD(H_0)$ is the domain of another
self-adjoint extension $H_n$ of the operator $\smo$. The extension $H_n$ is 
defined by the pair of defect subspaces $\gN_{i+nb}$ and $\gN_{-i+nb}$, and by 
the unitary operator $V_0^{(n)}$ in the space $\gN_{-i+nb}$. This operator is 
defined by the condition that its matrix with respect to the basis 
$\{U^n\psi_j\}$ coincides with the matrix  $V_0$. It is easily seen that 
$V_0^{(n)}=U^nV_0U^{*n}|\gN_{-i+nb}$.

The extension $H_n$ can be also  characterized in terms of the defect
 subspaces $\gN_i$ and $\gN_{-i}$ and the unitary operator $V_n$ acting on 
$\gN_{-i}$. In order to do it it is sufficiently to find the operator $V_n$
 from the system of equations
\begin{equation*}
\varphi_j-V_n\psi_j=f_{0,j}+\sum_{k}(U^n\varphi_k-V_0^{(n)}U^n\psi_k)
\alpha_{kj},\qquad f_{0,j}\in\gD(\smo), j=1,2\ldots m.
\end{equation*}
Let us introduce the following $m\times m$ matrices:
\begin{eqnarray}\label{matrices}
A_n=[(U^n\psi_k,\psi_l)]_{k,l=1}^m, \quad B_n=[(U^n\varphi_k,\psi_l)],\\ 
C_n=[(U^n\psi_k,\varphi_l)]_{k,l=1}^m,\quad
D_n=[(U^n\varphi_k,\varphi_l)]_{k,l=1}^m.
\end{eqnarray} 
Then the matrix of operator $V_n$ with respect to the basis $\{\psi_j\}$ is
defined by the expression
\begin{equation}\label{group}
V_n=T_n(V_0)=-\left[(nb-2i)A_nV_0-nbB_n\right]\left[nbC_nV_0-(nb+2i)D_n\right]^{-1}.
\end{equation}    
Putting $T_0=id-$ the identity mapping, we obtain the family 
$\Gamma={T_n,n\in\mZ}$ of mappings of the set $U(m)$ of $m\times m$ unitary 
matrices into itself. By its construction the mappings $T_n$ posses the 
property $T_n(T_m(\cdot))=T_{n+m}(\cdot)$. Therefore the family $\Gamma$ is a
group. 

From the Corollary 1 we obtain that if for some initial unitary matrix  $V_0$
the trajectory $\{T_k(V_0)\}_{k=-\infty}^{\infty}$ is periodic,that is
$T_n(V_0)=V_0$ for some positive integer $n$, than it is periodic for any
other initial matrix with the same period $n$. In such a situation the 
operator $\smo$ admits $(U,nb)-$ periodic self-adjoint extension, where $n$ is
the period of the trajectory of an initial unitary matrix $V_0$. We
reformulate this property as a property of the group $\Gamma$:

\begin{Prp} Let $\smo$ be a $(U,b)-$periodic prime symmetric operator with
index of defect $(m,m)$ and $\Gamma$ be the associated group of mappings of
the set $U(m)$ into itself, defined by (\ref{matrices}-\ref{group}). Then the operator $\smo$
admits periodic self-adjoint extension if and only if the group $\Gamma$ is
cyclic. 
\end{Prp}
\vskip .5truecm

{\bf Examples.} \vskip .5truecm
(a) Let $h(\lambda)$ be a nonnegative bounded function which has the period
$b$. Put $d\sigma(\lambda)=h(\lambda)/(1+\lambda^2)d\lambda$ and use the
definition (\ref{Wf}). Then the corresponding  function has the period
$b$. In particular, for $h(\lambda)=1+\sin{\lambda}$, 
\begin{equation*}
F(z)=i+e^{iz}-e^{-1}.
\end{equation*}
The function $F(z)$ has  the period $2\pi$. It is the Weyl-Titchmarsh function 
of the pair $(\smo, H)$ defined  by the formulas (\ref{domself}),
(\ref{self}),(\ref{ds}).
\vskip .5truecm
(b) Let $\gH=L^2_m[0,l]$, and the operator $\smo$ is defined as following:

Its domain is the set of all absolutely continuous functions 
$f(t)=\{f_k(t)\}_{k=1}^m\in\gH$, such that $f^{\prime}\in\gH$, $f(0)=f(l)=0$;
\begin{equation}\label{examplea}
\smo f(t)=i\frac{df}{dt}.
\end{equation} 
The operator $\smo$ has defect index $(m,m)$. The defect subspace $\gN_i$ is
generated by the columns of the $m\times m$ matrix 
$\exp{(t)}I_m$. There is one-to-one correspondence between set of self-adjoint
extensions of $\smo$ and $m\times m$ unitary matrices $V$.
Any self-adjoint extension $H_V$ of $\smo$ is obtained  as 
follows: 

Its domain is set of all absolutely continuous functions $f$ from
 $L^2_m[0,l]$,  such that $f^{\prime}\in L^2_m[0,1]$, and $f(0)=Vf(l)$, 
where $V$ is a unitary matrix in $\mC^m$. For the pair $(\smo, H_V)$ the 
Weyl-Titchmarsh function $M_{\smo, H_V}$ is equal to 
\begin{equation}
M_{\smo, H_V}(z)=-iI_m+\frac{2i}{e^{2l}-1}
(e^{l(1-iz)}-1)(I_m-e^{-izl}V)^{-1}(I_m-e^lV).
\end{equation}
This function has the period $2\pi/l$. Therefore, the operators 
(\ref{examplea}) and $H$ form a $2\pi/l-$periodic pair and the same is true
for
 any other self-adjoint extension of (\ref{examplea}). The unitary 
operator $U$, such that $U\smo U^*=\smo -(2\pi/l)I$, and similar equality for 
$H$ is the operator of multiplication by $\exp{(-2\pi it/l)}$.   
\vskip.3truecm
(c) More generally, consider the operator  $\smo_1=id/dt+V(t)$ on $L^2_m[0,1]$ 
with the same domain that above.  $V$ is a hermitian, bounded measurable
matrix function which satisfies condition $V(0)=V(l)$. Then the operator 
$\smo_1$ is symmetric with index of defect  $(m,m)$. Let $H_1$ be its 
self-adjoint extension. Then the Weyl-Titchmarsh function $M_{\smo_1, H_1}(z)$ 
has the period $2\pi/l$.
\vskip .3truecm
According to well-known theorem by M.Livsic \cite{Liv} 
a prime symmetric operator with index of defect (1,1) which admits a quasi-
hermitian extension $smo_v$ without spectrum in the finite complex plane is 
unitarily equivalent to operator, described in example (b) with $m=1$ for 
$l=2tr(\Im\smo_v^{-1})>0$.
Therefore, we have the following statement.
\begin{Thr} Let $\smo$ be a prime symmetric operator with index of defect
(1,1), and $H$ be a self-adjoint extension of $\smo$. Suppose that $\smo$
admits quasi-self-adjoint extension $\smo_v$ without spectrum. Then  the
Weyl-Titichmarsh function  $\WF(z)$ of the pair $(\smo,H)$ is a periodic 
one. Its period  is equal to $\pi/tr(\Im\smo_v^{-1})$.
\end{Thr}
This theorem does not admit generalization for the case of larger defect 
numbers. Indeed, let $\gH=L^2[0,l]$, and let $0<\xi<l$. Consider the symmetric
operator $\smo$ on $\gH$, defined as following:

\noindent The domain $\gD(\smo)$ is the set of all  functions $f(t)$ which 
are absolutely continuous for $0<t<\xi$ and $\xi<t<l$, $f^{\prime}\in\gH$, and 
$f(0)=f(\xi)=f(l)=0$. For $f\in\gD(\smo)$ $\smo f=idf/dt$.
The index of defect of $\smo$ is equal $(2,2)$. This operator admits
quasi-self-adjoint extension $\smo_v$ without spectrum, and $\smo_v^{-1}$ is
dissipative and unicellular \cite{Br}.
The operator $\smo$ is 
isomorphic to the direct sum $\smo_1\bigoplus\smo_2$ of two first order 
differential operators with zero boundary conditions on $[0,\xi]$ and 
$[\xi,l]$ respectively. Let $H$ be the self-adjoint extension of
$\smo_1\bigoplus\smo_2$ obtained by imposing the following conditions: 
$f(0)=\omega_1f(\xi-0)$, $f(\xi+0)=\omega_2f(l)$, where
$|\omega_1|=|\omega_2|=1$. The the Weyl-Titchmarsh function $\WF(z)$ of the
pair $(\smo, H)$ is a $2\times 2$ diagonal matrix 
\begin{equation*}
\WF(z)=\begin{bmatrix}M_1(z)&0\\
                      0&M_2(z)
\end{bmatrix},
\end{equation*}
where
\begin{eqnarray*}
M_1(z)=-i+2i(e^{\xi(1-iz)}-1)(1-\omega_1e^{\xi})/[(e^{2\xi}-1)(1-\omega_1e^{-iz\xi})],
\\
M_2(z)=-i+2i(\omega_2e^l-e^{\xi})(e^le^{-(l-\xi)iz}-e^{\xi})/[(e^{2l}-e^{2\xi})(\omega_2-e^{-iz(l-\xi)})].
\end{eqnarray*}
$M_1$ has the period $2\pi/\xi$, function $M_2$ has the period
$2\pi/(l-\xi)$. Therefore, if $\xi/(l-\xi)$ is an irrational number, the
function $\WF$ is not a periodic.

\section{Operators With Constant Weyl-Titchmarsh Function.}
Let $H$ be a self-adjoint operator, and let $W(t), t\in\mR$ be the 
one -parametric group of unitary operators generated by $H$ 
($W(t)=\exp{(iHt)}$). If $H$ is a $(U,b)$-periodic operator, then the
following  commutative relation is fulfilled:
\begin{equation}\label{onegroup} 
UW(t)=e^{-itb}W(t)U.
\end{equation}
\vskip .3truecm
So far we have considered the Weyl-Titchmarsh functions, which are invariant 
under some fixed shift $b$ of the argument. Let $F(z)$ be a function whose
values are operators on $m-$dimensional space $\gN$,   which 
admits representation  (\ref{rfunction}) and invariant under arbitrary real 
shift, that is $F(z+s)=F(z)$ for any real $s$. In such a situation the
function  $F(z)$ is, of course, constant  in each half-plane, 
\begin{equation}\label{lebesgue}
F(z)=\begin{cases}iI_{\gN}&\text{  $z\in\mC_+$},\\
  -iI_{\gN}&\text{  $z\in\mC_-$}. 
\end{cases}
\end{equation}
These properties are fulfilled
if and only if $d\tau(\lambda)=\pi^{-1}d\lambda I_{\gN}$. 
  
We have $F(z)=M_{\tilde\smo,\tilde H}(z)$ for the  pair $(\tilde\smo,\tilde H)$
acting in the Hilbert space $\tilde\gH=L^2(\mR,\gN,\pi^{-1}d\lambda)$, where
\begin{eqnarray}\label{const1} 
&\gD(\tilde H)=\{f\in L^2(\mR,\gN,\pi^{-1}d\lambda)|\int_{\smR}(1+\lambda^2)\Vert
f(\lambda)\Vert^2_{\gN}d\lambda<\infty\};&\\
\label{const2}
&(\tilde Hf)(\lambda)=\lambda f(\lambda);&\\
\label{const3}
&\gD(\tilde\smo)=\{f\in\gD(H)|\int_{\smR}f(\lambda)d\lambda=0\};&\\
\label{const4}
&(\tilde\smo f)(\lambda)=\lambda f(\lambda).&
\end{eqnarray} 
According to the Theorem 3 for any real number $s$ there is a unitary operator
$\tilde V(s)$ on $L^2(\mR,\gN,\pi^{-1}d\lambda)$  such that 
$\tilde V(s)\tilde H\tilde V^*(s)=\tilde H-sI$, and 
$\tilde V(s)\tilde\smo V^*(s)=\tilde\smo-sI$. The operators $\tilde V(s)$ act as following: 
$(\tilde V(s)f)(\lambda)=f(\lambda-s)$. Therefore,  the family $\{\tilde V(s)\}$ is 
strongly continuous unitary group. If $\tilde W(t)=\exp{(it\tilde H)}$, then
\begin{equation}\label{comrl}
\tilde V(s)\tilde W(t)=e^{-ist}\tilde W(t)\tilde V(s),
\end{equation}
which is the Weyl's form of the canonical commutative relation.

\begin{Thr}. Let $\smo$ be a prime symmetric operator  with index of defect
  $(m,m)$, $m<\infty$, $H\supset\smo$ be its self-adjoint extension, 
and let $W(t)(=\exp{(itH)})$ be the unitary group 
generated by $H$. Then the following conditions are equivalent
\begin{enumerate}
\item
There exists a unitary group $V(s)$ of s-type operators such that $V(s)W(t)=e^{-its}W(t)V(s)$;
\item
The Weyl-Titchmarsh function $\WF(z)=iI_{\gN_i}$ for $z\in\mC_+$, and\\ 
$\WF(z)=-iI_{\gN_i}$ for $z\in\mC_-$, where $\gN_i$, $dim\gN_i=m$, is the 
defect subspace of $\smo$. 
\end{enumerate}
\end{Thr} 
Let $G$ be the self-adjoint operator such that $V(s)=\exp{(isG)}$. Then
condition {\it 1} means that
\begin{equation*}
[G,H]=iI
\end{equation*}
on a dense subset of $\gH$.

{\bf Proof.} We have proved that from the statement {\it 2} follows the
statement {\it 1}. Let the statement {\it 1} is fulfilled. Then for
$f\in\gD(H)$ it follows that \\$V(s)f\in\gD(H)$ for any $s\in\mR$, and 
$V(s)Hf=(H-sI)V(s)f$. It is not hard to show that last  condition along with 
and the assumption about special structure of operators $V(s)$ implies that
the operator $\smo$ is also $(U,s)-$periodic for any real $s$.
 Therefore the Weyl-Titchmarsh function of the 
pair $(\smo,H)$ is constant in upper half-plane and in lower half-plane.

\noindent The Theorem is proved.

The pair $(\smo, H)$ is unitarily equivalent to its functional model
$(\tilde\smo,\tilde H)$ given
by the formulas (\ref{const1}-\ref{const4}). In such representation the group 
$\tilde V(s)$, as it was pointed \\out above, can be selected as group of shifts, 
$(\tilde V(s)f)(\lambda)=f(\lambda-s)$.

Consider the case $m=1$. The group $\tilde W(t)(=\exp{(i\tilde Ht)})$ is 
 the group  of multiplication by $\exp{(i\lambda t)}$ in the space
$\gH=L^2(\mR,\pi^{-1}d\lambda)$,
and $(\tilde V(s)f)(\lambda)=f(\lambda-s)$. This statement follows form the fact that for
each $s$ the operator $V(s)$ satisfies $\tilde V(s)H=(\tilde H-sI)\tilde V(s)$, Proposition 1, and
the group property ($\tilde V(s_1+s_2)=\tilde V(s_1)\tilde V(s_2)$). Therefore,
we obtained the statement of the Stone-von Neumann theorem for degree of
freedom 1 (\cite {St}). 

Let $D$ be the selfadjoint operator, such that $\tilde V(s)=\exp{(iDs)}$. Then 
\begin{eqnarray}\label{diff1}
&\gD(D)=\{f\in L^2(\mR, \pi^{-1}d\lambda)|f\in AC(-\infty, \infty);
f^{\prime}\in L^2(\mR, \pi^{-1}d\lambda)\},&\\
\label{diff2}
&(Df)(\lambda)=if^{\prime}(\lambda).&
\end{eqnarray}
The operator $D$ is the selfadjoint extension of the operator $\mathcal{D}$
defined as 
\begin{multline}\label{diff3}
\gD(\mathcal D)=\{f\in L^2(\mR, \pi^{-1}d\lambda)|f\in AC(-\infty, \infty);\\
f^{\prime}\in L^2(\mR, \pi^{-1}d\lambda); f(0)=0\},
\end{multline}
\begin{equation}\label{diff4}
({\mathcal D}f)(\lambda)=if^{\prime}(\lambda).
\end{equation}

Again applying the Theorem 5, we obtain that the Weyl-Titchmarsh function of
the pair $(\mathcal D,D)$ is constant (this fact can be checked, of course, by
direct calculations.). If $D_{\omega}$ and $\tilde H_{\theta}$ be arbitrary
selfadjoint extensions of $\mathcal D$ and $\tilde\smo$ respectively, then,
according to the Corollary 1, the Weyl-Titchmarsh functions $M_{\tilde\smo,
  \tilde H_{\theta}}(z)$ and $M_{\mathcal D, D_{\omega}}(z)$ are constant.
Therefore pair $(\tilde\smo, \tilde H_{\theta})$ is unitarily equivalent to the pair
$(\tilde\smo, \tilde H)$, and pair $(\mathcal D, D_{\omega})$ is unitarily equivalent to
the pair $(\mathcal D, D)$.

We have 
\begin{equation}\label{ht1}
\gD(\tilde H_{\theta})=\{f|f(\lambda)=f_0+(\frac{1}{\lambda-i}-\frac{\theta}{\lambda+i})z\},
\end{equation}
where $f_0\in\gD(\smo)$, $|\theta|=1$, and $z\in\mC$.
\begin{equation}\label{ht2}
(\tilde H_{\theta}f)(\lambda)=\lambda
f_0(\lambda)+i[1/(\lambda-i)+\theta/(\lambda+i)]z, 
\end{equation}
and $\tilde H=\tilde H_1$.
  The unitary operator $\Gamma_{\theta}$ such that 
$\tilde H_{\theta}=\Gamma_{\theta}\tilde H_1\Gamma_{\theta}^*$ acts as following: 
$(\Gamma_{\theta}f)(\lambda)=\theta\hat f_+(\lambda)+\hat f_-(\lambda)$, where
$f=\hat f_++\hat f_-$ is the (unique) representation of  function $f\in
L^2(\mR,d\lambda)$ as the sum of functions  $\hat f_+\in H^2_+$ and $\hat
f_-\in H^2_-$. Since $1/(\lambda-i)\in H^2_-$, and $1/(\lambda+i)\in H^2_+$, we
 need to show that $\Gamma_{\theta}\gD(\smo)\subset \gD(\smo)$. For $f\in
\gD(\smo)$ we have 
\begin{equation*}
f(\lambda)=\frac{1}{\sqrt{2\pi}}\int\limits_{-\infty}^{\infty}e^{i\lambda
  t}F(t)dt,
\end{equation*}
where $F\in L^2(\mR,dt)$, $F^{\prime}\in L^2(\mR,dt)$, and $F(0)=0$. 
\begin{equation}\label{gop}
(\Gamma_{\theta}f)(\lambda)\frac{1}{\sqrt{2\pi}}\int\limits_{-\infty}^{\infty}e^{i\lambda
  t}F(t)[\theta\chi_+(t)+\chi_-(t)]dt, 
\end{equation}
where $\chi_{\pm}$ are indicators functions of the positive and negative semaxes
respectively. The integrand of the last expression 
is equal to zero at $t=0$, therefore $\Gamma_{\theta}f\in\gD(\smo)$. It is
also clear that $\widehat{(\smo f)}_{\pm}=\lambda\hat f_{\pm}$, and
$\Gamma_{\theta}^*=\Gamma_{\bar\theta}$.

For the operator $D_{\omega}$ we have
\begin{multline}\label{dt1}
\gD(D_{\omega})=\{f\in L^2(\mR, d\lambda)|f\in AC([-R,0])\cap AC([0,R])\forall
R>0; \\
f(0_-)=\omega f(0_+),|\omega|=1; f^{\prime}\in L^2(\mR,d\lambda)\}
\end{multline}
\begin{equation}\label{dt2}
(D_{\omega}f)(\lambda)=if^{\prime}(\lambda),
\end{equation}
and $D=D_1$. 

The unirtary operator $J_{\omega}$ such that
$D_{\omega}=J_{\omega}D_1J_{\omega}^*$ acts as following:
\begin{equation}\label{jop}
(J_{\omega}f)(\lambda)=[\chi_-(\lambda)+ \omega\chi_+(\lambda)]f(\lambda),
\end{equation}
$J_{\omega}^*=J_{\bar\omega}$. 

From (\ref{gop}) and (\ref{jop}) it follows that
$\Gamma_{\theta}J_{\omega}=J_{\omega}\Gamma_{\theta}$.

Let $\tilde W_{\theta}$ be the unitary group generated by $\tilde H_{\theta}$, and 
$\tilde V_{\omega}(s)$ be the unitary group generated by $D_{\omega}$. It is not hard
to describe their actions. For example, the group $\tilde V_{\omega}(s)$ acts as following: 
  
for $s>0$
\begin{eqnarray*}
(\tilde V_{\omega}(s)f)(\lambda)=
\begin{cases}f_-(\lambda-s)&\text{$\lambda<0$}\\
\omega f_-(\lambda-s)&\text{$0\le\lambda\le s$}\\
f_+(\lambda-s)&\text{$\lambda\ge s$}
\end{cases}
\end{eqnarray*}
and for $s<0$
\begin{eqnarray*}
(\tilde V_{\omega}(s)f)(\lambda)=
\begin{cases}f_-(\lambda-s)&\text{$\lambda<s$}\\
\bar\omega f_+(\lambda-s)&\text{$s\le\lambda<0$}\\
f_+(\lambda-s)&\text{$\lambda\ge 0$}
\end{cases}
\end{eqnarray*}
    
It is clear, that $\Gamma_{\theta}D_1=D_1\Gamma_{\theta}$, and
$J_{\omega}H_1=H_1J_{\omega}$.

\begin{Prp} Let $\tilde H_{\theta}$ and $D_{\omega}$ be the operators defined by 
(\ref{ht1}-\ref{ht2}) and (\ref{dt1}-\ref{dt2}) respectively. Then for the
unitary groups $\tilde W_{\theta}(t)$ and $\tilde V_{\omega}(s)$ generated by
$\tilde H_{\theta}$ and $D_{\omega}$ respectively the H. Weyl commutative relation 
(\ref{comrl}) is fulfilled, that is 
\[\tilde V_{\omega}(s)\tilde W_{\theta}(t)=e^{-its}\tilde W_{\theta}(t)\tilde V_{\omega}(s)\]
\end{Prp}
The proposition follows from the following chain of equalities where above
mentioned properties of the operators $\Gamma_{\theta}$, $J_{\omega}$, $D_1$,
and $\tilde H_1$ are used:
\begin{multline*}
\tilde V_{\omega}(s)\tilde W_{\theta}(t)=J_{\omega}\tilde
V_1(s)J^*_{\omega}\Gamma_{\theta}\tilde W_1(t)\Gamma_{\theta}^*=
J_{\omega}\Gamma_{\theta}\tilde V_1(s)\tilde W_1(t)\Gamma_{\theta}^*J_{\omega}^*=\\
e^{-ist}J_{\omega}\Gamma_{\theta}\tilde W_1(t)\tilde V_1(s)\Gamma_{\theta}^*J_{\omega}^*=
e^{-ist}\Gamma_{\theta}\tilde W_1(t)\Gamma_{\theta}^*J_{\omega}\tilde V_1(s)J_{\omega}^*\\
=e^{-ist}\tilde W_{\theta}(t)\tilde V_{\omega}(s).
\end{multline*}

Last proposition admits reformulation in abstract form.
\begin{Prp}. Let $F_1$ and $G_1$ be self-adjoint operators with simple
  spectra acting in a Hilbert space $\gH$,   and
  corresponding unitary groups $ V_1(s)(=\exp{(iF_1s)})$ and $ W_1(t)(=\exp{(iG_1t)})$ 
satisfy (\ref{comrl}). Then: 
\begin{enumerate}
\item
There are prime symmetric operators $F_0$ and $G_0$ which have index of
defect $(1,1)$ such that $F_0\subset F_1$ and $G_0\subset G_1$;
\item
For any other self-adjoint extensions $F_{\omega}$
  and $G_{\theta}$ of the operators $F_0$ and $G_0$ respectively the
  corresponding unitary groups $ V_{\omega}(s)$ and $ W_{\theta}(t)$ also satisfy
 (\ref{comrl});
\item
There exists the unitary operator $U_{\theta\omega}:\gH\to
L^2(\mR,\pi^{-1}d\lambda)$ such that
$F_{\omega}=U_{\theta\omega}^*D_{\omega}U_{\theta\omega}$,
$G_{\theta}=U_{\theta\omega}^*\tilde H_{\theta}U_{\theta\omega}$,
$F_0=U_{\theta\omega}^*{\mathcal D}U_{\theta\omega}$, and 
$G_0=U_{\theta\omega}^*\tilde\smo U_{\theta\omega}$. 
\end{enumerate}
\end{Prp}
This proposition follows from the Stone-Von Neumann Theorem and previous
consideration. It also gives some refinement of the Stone-von Neumann's Theorem. 
The case $\omega=\theta=1$ is the most well-known. It corresponds to the
operators of momentum and coordinate in quantum mechanichs.
\vskip .3truecm 
Consider one more example of the pair  with constant
Weyl-Titchmarsh function.  Let $\gH=L^2(\mR, dt)$ and the self-adjoint operator
 is defined by the differntial expression 
\begin{equation}\label{shrd}
Lf=-\frac{1}{\gamma}\frac{d^2f}{dx^2}+xf,
\end{equation}
where $\gamma$ is a real constant. Corresponding self-adjoint operator
describes the particle in uniform electrical field.
This operator via Fourier transform is unitarily equivalent to the
self-adjoint operator $H$ defined as
\begin{multline*}
(Hf)(t)=i\frac{df}{dt}+\frac{1}{\gamma}t^2f(t);\\
\gD(H)=\{f\in L^2(\mR, dt)|f\in AC(-\infty,\infty),f^{\prime}\in
L^2(\mR,dt),\\t^2f(t)\in L^2(\mR,dt)\}.
\end{multline*}
Define the operator $\smo$ as following
\begin{multline*}
\gD(\smo)=\{f\in L^2(\mR, dt)|f\in AC(-\infty,0]\cup[0,\infty),f(0)=0,f^{\prime}\in
L^2(\mR,dt),\\t^2f(t)\in L^2(\mR,dt)\};
\end{multline*}
\begin{equation*}
(\smo f)(t)=i\frac{df}{dt}+\frac{1}{\gamma}t^2f(t).
\end{equation*}
The operator $\smo$ is a symmetric operator with index of defect $(1,1)$, and
$H$ is the selfadjoint extension of $\smo$.  
For any real $s$ define a unitary operator $U_s$ on $\gH$  by 
$(U_sf)(t)=e^{ist}f(t)$. Then we have $U_s\gD(\smo)=\gD(\smo)$, 
$U_s\gD(H)=\gD(H)$, and $U_sHU_s^*=(H- sI)$
that is the pair  $(\smo,H)$ is $(U_s, s)$-periodic.  
From the Theorem 5 it follows now that the Weyl-Titchmarsh function of the 
pair $(\smo,H)$ is constant in each half-plane. Therefore, the operator 
$H$ is unitarily equivalent to the operator of multiplication in 
$L^2(\mR, dt)$. 

Thus, the self-adjoint operator, generated by the differential expression
(\ref{shrd}) and its appropriate symmetric restriction have the constant 
Weyl-Titchmarsh function.

Let $V$ be a bounded, measurable, periodic,  real valued periodic
function. Without loss of generality we assume that the period of $V$ is $2\pi$.
The  Fourier series of $V$
\begin{equation*}
\sum\limits_{k=-\infty}^{\infty}\hat V(k)e^{ikx}
\end{equation*}
converges to $V(x)$ a.e., where $\hat V(k)$ are the Fourier coefficients of
the function $V$.

Consider the self-adjoint operator 
\begin{equation*}
L_1=L+V.
\end{equation*}
Again Fourier transform gives that the operator $L_1$ is unitarily equivalent
to the operator
\begin{equation*}
H_1f=i\frac{df}{dt}+\frac{1}{\gamma}t^2f+\sum\limits_k\hat V(k)f(t+k).
\end{equation*}
Operator $H_1$ is the selfadjoint extension of the symmetric operator
$\smo_1$with the same domain that the operator $\smo$ above. Now 
we have 
\begin{equation*}
U_sH_1f-H_1U_sf=-se^{ist}f+e^{ist}\sum\limits_{k}\hat V(k)(1-e^{isk})f(t+k),
\end{equation*}
and similar expression for $U_s\smo_1-\smo_1U_s$. 
Putting $s=2\pi$, we see that \\$U_{2\pi}H_1-H_1U_{2\pi}=-2\pi U_{2\pi}$, and
similar equation for $\smo_1$. Therefore, the pair $(\smo_1,H_1)$ is $2\pi-$
periodic. Thus the pair $(\mathcal L_1,L_1)$ where $\mathcal L_1$ is the symmetric
restriction of the Shr\"odinger operator $L_1$ with index of defect $(1,1)$
(inverse Fourier Transform of $\smo_1$) has the $2\pi-$periodic Weyl-Titchmarsh function.

%If $V(x)=\sin{bx}$, the the Weyl-Titichmarsh function $M_{\smo_1,H_1}$ has the
%period $2\pi/b$.
\vskip .5truecm
{\bf Acknowledgment.} Authors are very grateful to Fritz Gesztesy and
Konstantin A. Makarov for their suppot and stimulating discussions.
% The Appendices part is started with the command \appendix;
% appendix sections are then done as normal sections
% \appendix

% \section{}
% \label{}

% Bibliographic references with the natbib package:
% Parenthetical: \citep{Bai92} produces (Bailyn 1992).
% Textual: \citet{Bai95} produces Bailyn et al. (1995).
% An affix and part of a reference:
%   \citep[e.g.][Ch. 2]{Bar76}
%   produces (e.g. Barnes et al. 1976, Ch. 2).

\end{document}